\documentclass{amsart}
\usepackage[utf8]{inputenc}
\usepackage{bussproofs}
\usepackage{stmaryrd}
\usepackage{amsmath}
\usepackage{amsthm}
\usepackage{amssymb}
\usepackage{hyperref} 

\theoremstyle{plain}

\newtheorem*{normthm}{Normalization Theorem}
\newtheorem*{soundthm}{Soundness Theorem}

\newtheorem{proposition}{Proposition}
\newtheorem{corollary}{Corollary}

\theoremstyle{definition}
\newtheorem{definition}{Definition}

\theoremstyle{remark}
\newtheorem{remark}{Remark}

\usepackage{float}
\floatstyle{boxed}
\restylefloat{table}
\restylefloat{figure}

\usepackage{array}
\newcolumntype{L}[1]{>{\raggedright\let\newline\\\arraybackslash\hspace{0pt}}m{#1}}

\DeclareMathOperator{\KleeneT}{T}
\DeclareMathOperator{\RecC}{R}

\newcommand{\CTzero}{CT$_0$}

\newcommand{\AC}[2]{AC$^{{#1}{#2}}$}

\newcommand{\sr}[2]{{{#1} \text{ \textbf{mr} } {#2}}}
\newcommand{\srt}[2]{{{#1} \text{ \textbf{mrt} } {#2}}}
\newcommand{\HAomega}{HA$^\omega$}

\newcommand{\bn}{\mathbb{N}}

\newcommand{\Tplus}{T$^+$}

\newcommand{\PAomegaplus}{HA$^{\omega+}$}

\DeclareMathOperator{\Proof}{Proof}

\DeclareMathOperator{\hyp}{\textsf{hyp}}
\DeclareMathOperator{\wkn}{\textsf{wkn}}
\DeclareMathOperator{\lam}{\textsf{lam}}
\DeclareMathOperator{\app}{\textsf{app}}

\DeclareMathOperator{\pair}{\textsf{pair}}
\DeclareMathOperator{\fst}{\textsf{fst}}
\DeclareMathOperator{\snd}{\textsf{snd}}
\DeclareMathOperator{\zero}{\textsf{zero}}
\DeclareMathOperator{\suc}{\textsf{succ}}
\DeclareMathOperator{\rec}{\textsf{rec}}

\DeclareMathOperator{\shift}{\textsf{shift}}

\newcommand{\Hyp}{\hyp}
\newcommand{\Wkn}[1]{\wkn{#1}}
\newcommand{\Lam}[1]{\lam{#1}}
\newcommand{\App}[2]{\app({#1},{#2})}

\newcommand{\Pair}[2]{\pair({#1},{#2})}
\newcommand{\Fst}[1]{\fst{#1}}
\newcommand{\Snd}[1]{\snd{#1}}

\newcommand{\Suc}[1]{\suc{#1}}
\newcommand{\Rec}[3]{\rec({#1},{#2},{#3})}

\newcommand{\Shift}[1]{\shift{#1}}

\newcommand{\srConj}[4]{(\sr{\reify{}{}{\eval{}{}{\Fst{#1}}{#4}}}{#2}) \wedge (\sr{\reify{}{}{\eval{}{}{\Snd{#1}}{#4}}}{#3})}
\newcommand{\srImpl}[5]{\forall {#2}([\sr{\reify{}{}{\eval{}{}{#2}{#5}}}{#3}] \to [\sr{\reify{}{}{\eval{}{}{\App{#1}{#2}}{#5}}}{#4}])}
\newcommand{\srtImpl}[5]{\forall {#2}([\srt{\reify{}{}{\eval{}{}{#2}{#5}}}{#3}] \to [\srt{\reify{}{}{\eval{}{}{\App{#1}{#2}}{#5}}}{#4}])}
\newcommand{\srAll}[5]{\forall {#2}^{#3}(\sr{\reify{}{}{\eval{}{}{\App{#1}{#2}}{#5}}}{#4})}
\newcommand{\srEx}[3]{\sr{\reify{}{}{\eval{}{}{\Snd{#1}}{#3}}}{{#2}(\reify{}{}{\eval{}{}{\Fst{#1}}{#3}})}}

\DeclareMathOperator{\FV}{FV}
\DeclareMathOperator{\toptt}{tt}

\newcommand{\axc}[1]{\AxiomC{#1}}
\newcommand{\uic}[2]{\RightLabel{\small{#2}}\UnaryInfC{#1}}
\newcommand{\bic}[2]{\RightLabel{\small{#2}}\BinaryInfC{#1}}
\newcommand{\tic}[2]{\RightLabel{\small{#2}}\TrinaryInfC{#1}}

\newcommand{\reify}[3]{^{#1\!\!}\downarrow^{#2}{#3}}
\newcommand{\reflect}[3]{^{#1\!\!}\uparrow^{#2}{#3}}
\newcommand{\eval}[4]{^{#1}\llbracket {#3}\rrbracket^{#2}_{#4}}
\newcommand{\evalrec}[4]{^{#1}\{ {#3}\}^{#2}_{#4}}

\DeclareMathOperator{\er}{\textsf{e}}
\newcommand{\evdash}{\vdash_{\!\!\!\text{e}}}
\newcommand{\rvdash}{\vdash_{\!\!\!\text{r}}}
\newcommand{\gerefl}{\ge_{\text{refl}}}
\newcommand{\gecons}{\ge_{\text{cons}}}
\newcommand{\rvdashge}[2]{\ulcorner{#2}\urcorner^{#1}}
\newcommand{\evdashge}[2]{\llcorner{#2}\lrcorner_{#1}}
\newcommand{\Vdashge}[2]{\lceil{#2}\rceil^{#1}}
\newcommand{\Vvdashge}[2]{\lceil\!\!\lceil{#2}\rceil\!\!\rceil^{#1}}
\newcommand{\sVdashge}[2]{\lfloor{#2}\rfloor_{#1}}
\newcommand{\return}[1]{\eta{#1}}
\newcommand{\run}[1]{\mu{#1}}
\newcommand{\sVdash}{\Vdash_{\!\!\!\text{s}}}

\DeclareMathOperator{\proj}{proj}

\title[Interpretation of Sigma-2 Analysis in System T]{An interpretation of the Sigma-2 fragment \\of classical Analysis in System T}

\author{Danko Ilik}

\begin{document}


\begin{abstract}
We show that it is possible to define a realizability interpretation for the $\Sigma_2$-fragment of classical Analysis using Gödel's System T only. This supplements a previous result of Schwichtenberg regarding bar recursion at types 0 and 1 by showing how to avoid using bar recursion altogether. Our result is proved via a conservative extension of System T with an operator for composable continuations from the theory of programming languages due to Danvy and Filinski. The fragment of Analysis is therefore essentially constructive, even in presence of the full Axiom of Choice schema: Weak Church's Rule holds of it in spite of the fact that it is strong enough to refute the formal arithmetical version of Church's Thesis.
\end{abstract}



\maketitle

\section{Introduction}	
\label{sec:intro}

In the middle of the 20\textsuperscript{th} century, Kurt Gödel showed how to give a computational interpretation and a relative consistency proof of intuitionistic Arithmetic via his System T of equations between functionals definable by primitive recursion (in higher types) \cite{Godel1941,Godel1958}. Thanks to the fact that the induction axiom intuitionistically proves its own double negation translation, the interpretation also applies to classical Arithmetic. However, since the Axiom of Choice,
\begin{equation}
  \label{eq:ac}
  \tag{AC}
  \forall x\exists y A(x,y) \to \exists f\forall x A(x, f(x)),
\end{equation}
does not intuitionistically prove its double negation translation, the interpretation does not apply to classical Analysis. Gödel was of course aware of this fact and suggested \cite[\S 2.43]{Kreisel1959} that an extension of System T is needed which can interpret the logical schema,
\[
\neg\neg\forall x(A(x) \vee \neg A(x)),
\]
known as Kuroda's Conjecture \cite{Kuroda1951}, intuitionistically equivalent to the nowadays better known schema of Double Negation Shift,
\begin{equation}
  \label{eq:dns}
  \tag{DNS}
  \forall x\neg\neg B(x) \to \neg\neg\forall x B(x).
\end{equation}
Gödel must have also been aware of the difficulty involved in giving a computational interpretation to Kuroda's Conjecture, for already his 1941 lecture at Yale \cite{Godel1941} considers the special case when the formula $A(x)$ is $\exists y \KleeneT(x,x,y)$ --- where $\KleeneT$ is Kleene's predicate verifying that the Turing machine with code $x$, when run on input $x$, terminates with code $y$ --- which directly proves (see \cite{Troelstra1973}) the \emph{negation} of the formal arithmetical version of Church's Thesis,
\begin{equation}
  \label{eq:ct0}
  \tag{\CTzero}
  \forall x^\bn\exists y^\bn A(x,y) \to \exists e^\bn\forall x^\bn\exists u^\bn (T(e,x,u) \wedge A(x,U(u))).
\end{equation}

In spite of that, Spector and Kreisel \cite{Spector1962,GodelPostscriptSpector} managed to give a computational interpretation of DNS by extending System T with bar recursion, the computational adequacy of which was shown using a formal version of Brouwer's principle of Bar Induction \cite{Spector1962,Howard1968}. This approach to extracting computational content from proofs in Analysis via an extension of the primitive recursive System T with a \emph{general recursive} schema, has been much refined over the 50 years since it appeared \cite{KreiselBR1976,BerardiBC1998,Kohlenbach2008,BergerBS2002,Berger2004,BergerO2005,BergerO2006,Seisenberger} and has been applied to obtain results in Analysis proper, notably in Kohlenbach's Proof Mining programme \cite{Kohlenbach2008}. 

Nevertheless, as Schwichtenberg showed already in 1979 \cite{Schwichtenberg1979}, higher type primitive recursion is closed over the schema of bar recursion at types 0 and 1, and since a previous analysis of Kreisel \cite[\S 12.2]{Spector1962} shows that those low types are sufficient for interpreting the classical Axiom of Choice for formulas of the form $\exists \alpha\in{\bn\to\bn}\forall x\in\bn A_0(\alpha,x)$, $A_0$-quantifier-free, we in fact known that we should need no more than the primitive recursive functionals themselves in order to give a realizability interpretation of the uniformly realizable part of Analysis. Yet, it has remained unclear up to this day how to avoid using general recursive schemata altogether.

One alternative is offered by the use of so called computational side-effects (control operators) from the theory of programming languages. Krivine \cite{Krivine2003} used a realizability interpretation based on a virtual machine which can execute lambda calculus terms, extended with a control operator and a special machine instruction called ``quote''. Herbelin \cite{Herbelin2012} gave a more direct approach based on rewrite rules for a type theory extended with a control operator and a coinductive treatment of the existential quantifier. Both approaches rely on a proper extension of System T that can give a computational interpretation to full classical logic. However, given that there are classically true arithmetic statements that do not have a recursive realizer, it is not clear what the meaning of control operators outside the $\Sigma^0_1$-fragment is.

In this paper, we show that computational side-effects are not needed in the language of realizers, that is, although they are conceptually essential, control operators can be seen as a \emph{meta-mathematical} technique. Proofs of the $\Sigma_2$-fragment of Analysis, with the full Axiom of Choice, are essentially constructive and realizable by System T terms only. The soundness of the interpretation relies on a form of Markov's Principle (the \textsc{Shift} rule), rather than full classical logic \cite{BerardiBC1998,Krivine2003,Herbelin2012} or continuity principles \cite{Spector1962,Howard1968,Berger2004,BergerO2005}.

\section{Conservative extension of System T with operators for composable continuations}
\label{sec:systemtplus}

The constructive interpretation of proofs of Section~\ref{sec:hacplus} is based on Gödel's System T in its lambda-calculus formulation \cite{Troelstra1973,SchwichtenbergW2012}. Nevertheless, we consider it a conceptual advantage to use an intermediate system, the System \Tplus\ obtained when System T is extended with a control operator for composable continuations, the so called call-by-name variant of the \textit{shift} operator of Danvy and Filinski \cite{DanvyF1990}. The control operator is a key conceptual ingredient that led us to the interpretation \cite{HerbelinMP,IlikThesis,Ilik2010,IlikN2014}.

The goal of this technical section is to prove the following conservativity theorem, as well as prove that suitable equations important for Section~\ref{sec:hacplus} hold of System~T\textsuperscript{(+)} (Proposition~\ref{prop:equalities}).

\begin{normthm} Every term of System \Tplus maps to a term in normal form of System T.
\end{normthm}

The \emph{types} of System~T\textsuperscript{(+)} are $\bn$, functions $\bn\to\bn$, functionals (ex. $(\bn\to\bn)\to\bn$), and Cartesian products of these. Types will be denoted by $\sigma, \tau$:
\[
\mathcal{T} \ni \sigma,\tau ::= \bn ~|~ \sigma\to\tau ~|~ \sigma*\tau.
\]

\emph{Terms} of System \Tplus\ are constants associated with a sequent $\sigma_1;\ldots;\sigma_n\vdash \tau$ (also written $\gamma\vdash \tau$ for $\gamma$ a finite ordered list of types) which means that a term is of type $\tau$ and the free variables that appear in the term are of types $\sigma_1,\ldots,\sigma_n$. The terms, marked in \textsf{sans-serif} face,  are defined inductively as follows.
\begin{gather*}
    \hyp  \frac{}{(\sigma ; \gamma) \vdash \sigma}\qquad
    \wkn  \frac{\gamma \vdash \sigma}{(\tau ; \gamma) \vdash \sigma}\qquad
    \lam  \frac{(\sigma ; \gamma) \vdash \tau}{\gamma \vdash \sigma \to \tau}\\ \\
    \app  \frac{\gamma \vdash \sigma \to \tau\quad\gamma \vdash \sigma}{\gamma \vdash \tau}\qquad
    \pair  \frac{\gamma \vdash \sigma\quad\gamma \vdash \tau}{\gamma \vdash \sigma * \tau}\qquad
    \fst  \frac{\gamma \vdash \sigma * \tau}{\gamma \vdash \sigma}\\ \\
    \snd  \frac{\gamma \vdash \sigma * \tau}{\gamma \vdash \tau}\qquad
    \zero  \frac{}{\gamma \vdash \bn}\qquad
    \suc  \frac{\gamma \vdash \bn}{\gamma \vdash \bn}\\ \\
    \rec  \frac{\gamma \vdash \bn\quad\gamma \vdash \sigma\quad\gamma \vdash \bn \to \sigma \to \sigma}{\gamma \vdash \sigma}\qquad
    \shift  \frac{(\bn \to \sigma \to \bn ; \gamma) \vdash \bn}{\gamma \vdash \sigma}
\end{gather*}
For example, the Ackermann function of type $\bn\to\bn\to\bn$ in lambda calculus notation,
\begin{align*}
  \mathrm{A} &:= \lambda m. \RecC m (\lambda n. n+1) (\lambda m'.\lambda u.\lambda n. \RecC n (u 1) (\lambda n'.\lambda w. u w)),
\end{align*}
where $\RecC$ is a constant such that
\begin{align*}
  \RecC 0 a b &= a\\
  \RecC (n+1) a b &= b n (\RecC n a b),
\end{align*}
is defined using the following term:
\begin{align*}
&   \lam\\
&\quad     (\rec \hyp (\lam (\suc \hyp))\\
&\quad\quad      (\lam\\
&\quad\quad\quad       (\lam\\
&\quad\quad\quad\quad        (\lam\\
&\quad\quad\quad\quad\quad         (\rec \hyp (\app (\wkn \hyp) (\suc \zero))\\
&\quad\quad\quad\quad\quad\quad          (\lam (\lam (\app (\wkn (\wkn (\wkn \hyp))) \hyp)))))))).
\end{align*}
This style of presenting lambda calculus terms, known as deBruijn convention, keeps the language of terms first-order, that is, it avoids problems related to handling variables as names that we have to keep track of externally to the system: 
$\hyp$ denotes the variable corresponding to the nearest preceding lambda abstraction $\lam$ or a $\shift$, $\Wkn{\Hyp}$ denotes the second most recent introduced variable, $\Wkn{\Wkn{\Hyp}}$ the third one, and so on. Seen as natural numbers ($\hyp$ is 0, $\wkn$ is successor), these are just so called deBruijn indices.

The terms of System T have a great computational potential, they can compute any higher-type primitive recursive function. There is, however, a subclass of terms which denote computations that have finished, the so called terms \emph{in normal form}. These will be sufficient for denoting programs and data extracted from proofs in classical Analysis in Section~\ref{sec:hacplus}. Technically, the normal forms are known as $\beta$-normal in reference to the associated $\beta$-reduction relation. Here we give a direct inductive characterization of normal terms ($\rvdash$),
\begin{gather*}
    \er  \frac{\gamma \evdash \sigma}{\gamma \rvdash \sigma}\qquad
    \lam  \frac{(\sigma ; \gamma) \rvdash \tau}{\gamma \rvdash \sigma \to \tau}\qquad
    \pair  \frac{\gamma \rvdash \sigma\quad\gamma \rvdash \tau}{\gamma \rvdash \sigma * \tau}\\
    \zero  \frac{}{\gamma \rvdash \bn}\qquad
    \suc  \frac{\gamma \rvdash \bn}{\gamma \rvdash \bn}
\end{gather*}
defined at the same time with the so called \emph{neutral} terms ($\evdash$),
\begin{gather*}
    \hyp  \frac{}{(\sigma ; \gamma) \evdash \sigma}\qquad
    \wkn  \frac{\gamma \rvdash \sigma}{(\tau ; \gamma) \evdash \sigma}\qquad
    \app  \frac{\gamma \evdash \sigma \to \tau\quad\gamma \rvdash \sigma}{\gamma \evdash \tau}\\ \\
    \fst  \frac{\gamma \evdash \sigma * \tau}{\gamma \evdash \sigma}\qquad
    \snd  \frac{\gamma \evdash \sigma * \tau}{\gamma \evdash \tau}\qquad
    \rec  \frac{\gamma \evdash \bn\quad\gamma \rvdash \sigma\quad\gamma \rvdash \bn \to \sigma \to \sigma}{\gamma \evdash \sigma}.
\end{gather*}
Neutral terms correspond to computation that are ``blocked'': neutral terms are those that contain open/free variables that block a $\beta$-reduction step from happening. 

A property that will be later used (proof of Corollary~\ref{cor:gamma1}) follows directly from the shape of normal forms: any closed normal term of type $\bn$, i.e. a term of type $\emptyset\rvdash \bn$, is actually a numeral, that is, built only from $\zero$ and $\suc$-terms. This follows because \emph{closed} normal terms cannot be neutral: a neutral term necessarily has at least one free variable.

We are now ready to state our first theorem precisely.
\begin{normthm} There is a normalization procedure $\reify{}{}{\eval{}{}{-}{}}$ such that, for every term $p$ of System \Tplus of type $\gamma\vdash\tau$, the term $\reify{}{}{\eval{}{}{p}{}}$ is a term in normal form of System T of the same type ($\gamma\rvdash\tau$).
\end{normthm}

This theorem is not a standard fact from the theory of lambda calculus with control operators. As a matter of fact, we present the first proof that control operators can be completely eliminated from System T.

Specialists will recognize the proof method as a \emph{normalization-by-evaluation} \cite{BergerS1991} or \emph{type-directed partial evaluation} \cite{Danvy1996} argument. Nevertheless, the addition of control operators requires us to perform the proof in so called \emph{continuation-passing style}; a similar technique has been used to provide a constructive completeness proof for non-minimal intuitionistic logic when disjunction and the existential quantifier are present \cite{IlikThesis,Ilik2013}.

The proof is constructive and can be formalized in a suitable predicative meta-theory such as Martin-Löf Type Theory (for a modern formulation see \cite{hottbook}). At her convenience, the reader may find this full mechanization of the proof in Agda notation in \cite{code}; one can also use this machine-checked constructive proof directly in order to compute, for example, that the term for the Ackermann function A(3,2) really evaluates to $\underbrace{\suc\cdots\suc}_{29 \text{ times }}\zero$.

\begin{proof}[Proof of Normalization Theorem]
Our goal is to define an \emph{evaluation} function,
\[
  \eval{}{}{(-)}{} : \gamma \vdash \sigma \Rightarrow \gamma \Vdash \sigma
\]
that maps a term $p$ of type $\gamma \vdash \sigma$ to a \emph{forcing}\footnote{The terminology ``forcing'' comes from similarity of our construction with Kripke models which does not intentionally refer to the forcing from Set Theory.} set $\gamma \Vdash \sigma$, together with a \emph{reification} function,
\[
\reify{}{}{(-)} : \gamma \Vdash \sigma \Rightarrow \gamma \rvdash \sigma,
\]
that extracts normal forms (without $\shift$!) from the forcing set. Composing the two function, $p \mapsto \reify{~}{}{\eval{}{}{p}{}}$, gives a normalization procedure. (The notation $a~ b~ c \mapsto d$ will be a compact form of $a\mapsto b\mapsto c\mapsto d$.)

We first need a precise way to speak about extensions of type contexts $\gamma$ (imposed by the fact that we want to normalize potentially open terms). This is formalized by the initial segment, or prefix, preorder $\ge$ as follows.
\begin{gather*}
  \gerefl\frac{}{\gamma\ge\gamma}\qquad
  \gecons\frac{\gamma_2\ge\gamma_1}{(\sigma;\gamma_2)\ge\gamma_1}
\end{gather*}
For example, the proof of $\tau_1;\tau_2;\gamma\ge \gamma$ will be denoted by $\gecons\gecons\gerefl$. The transitivity of $\ge$ is proven as the (right-associative) operation $(-)\cdot(-)$, defined by recursion on the construction of the proofs of $\gamma_3 \ge \gamma_2$ and $\gamma_2 \ge \gamma_1$.
\begin{align*}
  (-)\cdot(-) &: \gamma_3 \ge \gamma_2 \Rightarrow \gamma_2 \ge \gamma_1 \Rightarrow \gamma_3 \ge \gamma_1\\
  \ge_3 \cdot \gerefl &= \ge_3\\
  \gerefl \cdot \ge_2 &= \ge_2\\
  (\gecons \ge_3) \cdot \ge_2 &= \gecons (\ge_3 \cdot \ge_2)  
\end{align*}
In this definition, and henceforth, the notation $\ge_n$ will be used to denote a (hypothetical) proof of $\gamma_n\ge\gamma_k$. For example, given a proof $\ge_2$ of $\gamma_2\ge \sigma;\gamma_1$, one can prove $\gamma_2\ge\gamma_1$ by the denotation $\ge_2 \cdot \gecons\gerefl$.

We can now define precisely the \emph{forcing} set $\gamma\Vdash \sigma$. This is done simultaneously with the \emph{strong forcing} set $\gamma\sVdash \sigma$, itself defined inductively following the construction of the type $\sigma$.
\begin{align*}
  \gamma\Vdash \sigma &= \forall\gamma_1 \ge \gamma \left(\forall\gamma_2 \ge \gamma_1 \left(\gamma_2 \sVdash \sigma \Rightarrow \gamma_2 \rvdash \bn\right) \Rightarrow \gamma_1 \rvdash \bn\right)\\
  ~&~\\
  \gamma \sVdash \bn &= \gamma \rvdash \bn\\
  \gamma \sVdash (\sigma \to \tau) &= \forall \gamma' \ge \gamma (\gamma' \Vdash \sigma \Rightarrow \gamma' \Vdash \tau)\\
  \gamma \sVdash (\sigma * \tau) &= \gamma \Vdash \sigma \times \gamma \Vdash \tau\\
\end{align*}
We will also need $\gamma'\Vvdash\gamma$, the component-wise extension of the forcing relation defined as follows ($[]$ denotes the empty context).
\begin{align*}
  \gamma' \Vvdash [] &= \top\\
  \gamma' \Vvdash (\sigma ; \gamma) &= (\gamma' \Vdash \sigma) \times (\gamma' \Vvdash \gamma)
\end{align*}
The first equation defines the forcing of the empty context to be the singleton set $\top$, whose unique inhabitant is denoted $\toptt$. The symbol $\times$ constructs a Cartesian product i.e. pair type in the ambient type theory; the components of a Cartesian product can be accessed by the projection operations $\proj_1$ and $\proj_2$.

The following operations (lemmas) show that $\rvdash,\evdash,\Vdash,\sVdash,\Vvdash$ are all monotone with respect to the prefix preorder. The last two operations are defined by induction on types and contexts, respectively.
\begin{align*}
  \rvdashge{(-)}{(-)} &: \gamma_2 \ge \gamma_1 \Rightarrow \gamma_1 \rvdash \sigma \Rightarrow \gamma_2 \rvdash \sigma\\
  \rvdashge{\gerefl}{H} &= H\\
  \rvdashge{\gecons \ge_2}{H} &= \er (\wkn (\rvdashge{\ge_2}{H}))
\end{align*}
\begin{align*}
  \evdashge{(-)}{(-)} &: \gamma_2 \ge \gamma_1 \Rightarrow \gamma_1 \evdash \sigma \Rightarrow \gamma_2 \evdash \sigma\\
  \evdashge{\gerefl}{H} &= H\\
  \evdashge{\gecons \ge_2}{H} &= \wkn (\er (\evdashge{\ge_2}{H}))
\end{align*}
\begin{align*}
  \Vdashge{(-)}{(-)} &: \gamma_2 \ge \gamma_1 \Rightarrow \gamma_1 \Vdash \sigma \Rightarrow \gamma_2 \Vdash \sigma\\
  \Vdashge{\gerefl}{H} &= H\\
  \Vdashge{\gecons \ge_2}{H} &= \ge_3 \mapsto H (\ge_3 \cdot \gecons \ge_2)  
\end{align*}
\begin{align*}
  \sVdashge{(-)}{(-)}^\sigma &: \gamma_2 \ge \gamma_1 \Rightarrow \gamma_1 \sVdash \sigma \Rightarrow \gamma_2 \sVdash \sigma\\
  \sVdashge{\ge_2}{H}^{\sigma \to \tau}  &=  \ge_3 \mapsto H (\ge_3 \cdot \ge_2)\\
  \sVdashge{\ge_2}{(H_1 , H_2)}^{\sigma * \tau} &= \Vdashge{\ge_2}{H_1} , \Vdashge{\ge_2}{H_2}\\
  \sVdashge{\ge_2}{H}^{\bn} &= \rvdashge{\ge_2}{H}
\end{align*}
\begin{align*}
  \Vvdashge{(-)}{(-)}_{\gamma_2} &: \gamma_2 \ge \gamma_1 \Rightarrow \gamma_1 \Vvdash \gamma \Rightarrow \gamma_2 \Vvdash \gamma\\
  \Vvdashge{\ge_2}{H}_{[]} &= H\\
  \Vvdashge{\ge_2}{H}_{\sigma ; \gamma} &= \Vdashge{\ge_2}{\proj_1 H} , \Vvdashge{\ge_2}{\proj_2 H}_{\gamma}
\end{align*}

Finally, we also need lemmas relating the forcing sets and derivability, the \emph{return} ($\eta$) and \emph{run} ($\mu$) operations.
\begin{align*}
  \return{(-)} &: \gamma \sVdash \sigma \Rightarrow \gamma \Vdash \sigma & \run{(-)} &: \gamma \Vdash \bn \Rightarrow \gamma \sVdash \bn\\
  \return{H} &= \ge_1 \kappa \mapsto \kappa \gerefl \sVdashge{\ge_1}{H} & \run{H} &= H \gerefl (\ge_1 \alpha \mapsto \alpha) 
\end{align*}
Return shows that we can always lift a member of the strong forcing set to a member of the forcing set. Run shows that, whenever we have a member of a set forcing type $\bn$, we can actually obtain a term of System T in normal form from it; note that by definition the sets $\gamma\sVdash\bn$ and $\gamma\rvdash\bn$ are the same.

The \emph{reify} function $\reify{}{}{(-)}$ shows that we can actually run any forcing set, for any type $\sigma$, and not just for $\sigma = \bn$. It is defined by induction on the type.
\begin{align*}
  \reify{\gamma}{\sigma}{(-)} &: \gamma \Vdash \sigma \Rightarrow \gamma \rvdash \sigma\\
    \reify{\gamma}{\bn}{H} &= \run{H}\\
    \reify{\gamma}{\sigma \to \tau}{H} &= 
      \lam (\reify{\gamma}{\tau}{}\\
      &\qquad   (\ge_1 \kappa \mapsto\\
      &\qquad\qquad      H (\ge_1 \cdot \gecons \gerefl)\\
      &\qquad\qquad      (\ge_2 \phi \mapsto\\
      &\qquad\qquad\qquad         \phi \gerefl (\Vdashge{\ge_2 \cdot \ge_1}{\reflect{\sigma;\gamma}{\sigma}{hyp}}) \gerefl\\
      &\qquad\qquad\qquad         (\ge_3 \mapsto \kappa (\ge_3 \cdot \ge_2)))))\\
    \reify{\gamma}{\sigma * \tau}{H} &= 
      \pair\\
      &\qquad  \reify{\alpha}{\sigma}{}(\ge_1 \kappa \mapsto \\
      &\qquad\qquad\qquad H \ge_1 (\ge_2 \alpha \mapsto \proj_1 \alpha \gerefl (\ge_3 \mapsto \kappa (\ge_3 \cdot \ge_2))))\\
      &\qquad  \reify{\alpha}{\tau}{}(\ge_1 \kappa \mapsto \\
      &\qquad\qquad\qquad H \ge_1 (\ge_2 \alpha \mapsto \proj_2 \alpha \gerefl (\ge_3 \mapsto \kappa (\ge_3 \cdot \ge_2))))
\end{align*}
The case of arrow type forces us to define, at the same time with the reify function, the reflect function $\reflect{}{}{(-)}$:
\begin{align*}
    \reflect{\gamma}{\sigma}{(-)} &: \gamma \evdash \sigma \Rightarrow \gamma \Vdash \sigma\\
    \reflect{\gamma}{\bn}{p} &= \return{(\er p)}\\
    \reflect{\gamma}{\sigma \to \tau}{p} &= \return{
                          (\ge_2 \alpha \mapsto \reflect{\gamma}{\tau}{\App{\evdashge{\ge_2}{p}}{\reify{\gamma}{\sigma}{\alpha}}})}\\
    \reflect{\gamma}{\sigma * \tau}{p} &= \return{(\reflect{\gamma}{\sigma}{\fst p} , \reflect{\gamma}{\tau}{\snd p})}
\end{align*}
Note that this function needs as domain only the neutral terms.

Finally, we are ready to define the \emph{evaluation} function $\eval{}{}{-}{}$ that constructs a member of the forcing set for any input term of System \Tplus. The definition is by recursion on the construction of the term $p$.
\begin{align*}
  \eval{\gamma}{\sigma}{(-)}{(-)} &: \gamma \vdash \sigma \Rightarrow \forall\gamma' \Vvdash \gamma (\gamma' \Vdash \sigma)\\
  \eval{}{}{\hyp}{\rho} &= \proj_1 \rho\\
  \eval{}{}{\Wkn{p}}{\rho} &= \eval{}{}{p}{\proj_2 \rho}\\
  \eval{}{}{\Lam{p}}{\rho} &= \return{(\ge_1 \alpha \mapsto \eval{}{}{p}{(\alpha , \Vvdashge{\ge_1}{\rho})})}\\
  \eval{}{}{\App{p}{q}} \rho &= 
    \ge_1 \kappa \mapsto\\
  &\qquad      \eval{}{}{p}{\rho} \ge_1\\
  &\qquad      (\ge_2 \phi \mapsto\\
  &\qquad\qquad         \phi \gerefl (\eval{}{}{q}{\Vvdashge{\ge_2 \cdot \ge_1}{\rho}}) \gerefl\\
  &\qquad\qquad         (\ge_3 \mapsto \kappa (\ge_3 \cdot \ge_2)))\\
  \eval{}{}{\Pair{p}{q}}{\rho} &= \return{(\eval{}{}{p}{\rho} , \eval{}{}{q}{\rho})}\\
  \eval{}{}{\Fst{p}} \rho &= 
    \ge_1 \kappa \mapsto
        \eval{}{}{p}{\rho} \ge_1 (\ge_2 \alpha \mapsto \proj_1 \alpha \gerefl (\ge_3 \mapsto \kappa (\ge_3 \cdot \ge_2)))\\
  \eval{}{}{\Snd{p}}{\rho} &= 
    \ge_1 \kappa \mapsto
        \eval{}{}{p}{\rho} \ge_1 (\ge_2 \alpha \mapsto \proj_2 \alpha \gerefl (\ge_3 \mapsto \kappa (\ge_3 \cdot \ge_2)))\\
  \eval{}{}{\Shift{p}}{\rho} &= 
    \ge_1 \kappa \mapsto
        \run{}
        \eval{}{}{p}{
         \return{}
          (\ge_2 \nu \mapsto
             \return{
             (\ge_3 \alpha \mapsto \return{(\alpha \gerefl (\ge_4 \mapsto \kappa (\ge_4 \cdot \ge_3 \cdot \ge_2)))})})
          , \Vvdashge{\ge_1}{\rho}}\\
  \eval{}{}{\zero}{\rho} &= \return{(\zero)}\\
  \eval{}{}{\Suc{p}}{\rho} &= \return{(\Suc{(\run{\eval{}{}{p}{\rho}})})}\\
  \eval{}{}{\Rec{n}{a}{f}}{\rho} &= \ge_1\kappa\mapsto\eval{}{}{n}{\rho} \ge_1 (\ge_2 \nu \mapsto \evalrec{}{}{\nu}{\Vvdashge{\ge_2\cdot\ge_1}{\rho}} \gerefl (\ge_3 \mapsto \kappa (\ge_3\cdot\ge_2)))\\
  \text{ where } & \\
  \evalrec{}{}{\zero}{\rho'} &= \eval{}{}{a}{\rho'}\\
  \evalrec{}{}{\Suc{r}}{\rho'} &= \ge_1\kappa\mapsto \eval{}{}{f}{\Vvdashge{\ge_1}{\rho'}} \gerefl (\ge_2 \gamma\mapsto \gamma \gerefl (\eta\rvdashge{\ge_2\cdot\ge_1}{r}) \gerefl\\
  & \qquad (\ge_3 \delta \mapsto \evalrec{}{}{r}{\rho'} (\ge_3\cdot\ge_2\cdot\ge_1) (\ge_4\alpha\mapsto \delta \ge_4 (\eta\alpha) \gerefl \\
  & \qquad\qquad (\ge_5\mapsto \kappa (\ge_5\cdot\ge_4\cdot\ge_3\cdot\ge_2)))))\\
  \evalrec{}{}{\er{e}}{\rho'} &= \reflect{~}{}{\Rec{e}{\reify{}{}{\eval{}{}{a}{\rho'}}}{\reify{}{}{\eval{}{}{f}{\rho'}}}}
\end{align*}
Note that the argument $\rho$ is of type $\gamma' \Vvdash \gamma$.

For $\gamma'=\gamma$, such a $\rho$ can always be constructed by reflecting the term $\hyp$:
\begin{align*}
  \Uparrow{\gamma} &: \gamma \Vvdash \gamma\\
  \Uparrow{[]} &= tt\\
  \Uparrow{\sigma ; \gamma} &= (\reflect{}{\sigma}{\hyp}) , \Vvdashge{\gecons \gerefl}{\Uparrow{\gamma}}
\end{align*}

We have therefore shown that, given $p : \gamma\vdash \sigma$ of System \Tplus, there is a term in normal form $\reify{}{}{\eval{}{}{p}{\rho}} : \gamma\rvdash \sigma$ of System T, for every $\rho$, and in particular one such term is $\reify{}{}{\eval{}{}{p}{\Uparrow{\gamma}}}$.
\end{proof}

The following proposition characterizes the equational theory generated by the normalization procedure. It will be used in the proof of Soundness Theorem of Section~\ref{sec:hacplus}. This has also been machine checked and is available in Agda notation from \cite{code}.

\begin{proposition}\label{prop:equalities} The following definitional equalities hold,
  \begin{align}
    \label{eq:wkn}\reify{}{}{\eval{}{}{\Wkn{p}}{\alpha,\rho}} &= \reify{~}{}{\eval{}{}{p}{\rho}} &\text{ for } \alpha\in\tau;\gamma\Vdash\tau\\
    \label{eq:hyp}\reify{}{}{\eval{}{}{\hyp}{\alpha,\rho}} &= \reify{~}{}{\alpha} &\text{ for } \alpha\in \gamma\Vdash\tau\\
    \label{eq:fst}\reify{}{}{\eval{}{}{\Fst{\Pair{p}{q}}}{\rho}} &= \reify{~}{}{\eval{}{}{p}{\rho}} & \\
    \label{eq:snd}\reify{}{}{\eval{}{}{\Snd{\Pair{p}{q}}}{\rho}} &= \reify{~}{}{\eval{}{}{q}{\rho}} & \\
    \label{eq:lam}\reify{}{}{\eval{}{}{\App{\Lam{p}}{q}}{\rho}} &= \reify{~}{}{\eval{}{}{p}{\eval{}{}{q}{\rho},\rho}} & \\
    \label{eq:zero}\reify{}{}{\eval{}{}{\Rec{\zero}{p}{q}}{\rho}} &= \reify{~}{}{\eval{}{}{p}{\rho}} & \\
    \label{eq:suc}\reify{}{}{\eval{}{}{\Rec{\Suc{r}}{p}{q}}{\rho}} &= \reify{~}{}{\eval{}{}{\App{\App{q}{r}}{\Rec{r}{p}{q}}}{\rho}} & \\
    \label{eq:shift:N}\reify{}{\bn}{\eval{}{}{\Shift{p}}{\rho}} &= \reify{~}{\bn}{\eval{}{}{p}{\phi,\rho}} & \\
    \label{eq:hyp:cont}\reify{}{\bn}{\eval{}{}{\App{\App{\hyp}{x}}{y}}{\phi,\rho}} &= \reify{~}{\bn}{\eval{}{}{y}{\phi,\rho}} & 
  \end{align}
  where for the last two equations,
  \[\phi := \eta(\ge_2\nu\mapsto\eta(\ge_3\alpha\mapsto\eta(\mu\alpha))),\]
  and $x, y : \bn\to\bn\to\bn \vdash \bn$.
\end{proposition}
\begin{proof}
Equations (\ref{eq:wkn})--(\ref{eq:suc}) follow from the ones that hold already of the $\eval{}{}{-}{(-)}$ function. This is because, as an argument to the $\reify{}{}{(-)}$ function, the evaluation function is always applied to some $\ge_1$ and $\kappa$.
  \begin{align*}
    \eval{}{}{p}{\rho} \ge_1 \kappa &= \eval{~}{}{q}{\rho} \ge_1 \kappa\\
    {\eval{}{}{\Wkn{p}}{\alpha,\rho}} \ge_1 \kappa &= {\eval{}{}{p}{\rho}} \ge_1 \kappa\\
    {\eval{}{}{\hyp}{\alpha,\rho}} \ge_1 \kappa &= {\alpha} \ge_1 \kappa \\
    {\eval{}{}{\Fst{\Pair{p}{q}}}{\rho}} \ge_1 \kappa &= {\eval{}{}{p}{\rho}} \ge_1 \kappa \\
    {\eval{}{}{\Snd{\Pair{p}{q}}}{\rho}} \ge_1 \kappa &= {\eval{}{}{q}{\rho}} \ge_1 \kappa \\
    {\eval{}{}{\App{\Lam{p}}{q}}{\rho}} \ge_1 \kappa &= {\eval{}{}{p}{\eval{}{}{q}{\rho},\rho}} \ge_1 \kappa \\
    {\eval{}{}{\Rec{\zero}{p}{q}}{\rho}} \ge_1 \kappa &= {\eval{}{}{p}{\rho}} \ge_1 \kappa \\
    {\eval{}{}{\Rec{\Suc{r}}{p}{q}}{\rho}} \ge_1 \kappa &= {\eval{}{}{\App{\App{q}{r}}{\Rec{r}{p}{q}}}{\rho}} \ge_1 \kappa
  \end{align*}
These equations come out by unfolding the definition and occasionally using an $\eta$-equality step of the form $(\alpha \mapsto \phi \alpha) = \phi$.

Equations (\ref{eq:shift:N})--(\ref{eq:hyp:cont}) also follow by definition, this time reification being applied for only one concrete type, $\bn$.
\end{proof}

\section{A modified realizability interpretation of Analysis}
\label{sec:hacplus}

By a logical theory sufficient to formalize proofs of Analysis we have in mind the System \PAomegaplus\ of Figure~\ref{fig:haplus} together with the full axiom of choice schema $\text{\AC{}{}}=\cup_{\sigma,\tau\in\mathcal{T}}\text{\AC{\sigma}{\tau}}$,
\begin{equation}
  \label{eq:acsigmatau}
  \tag{\AC{\sigma}{\tau}}
  \forall x^\sigma\exists^\tau y A(x,y) \to \exists^{\sigma\to\tau} f\forall x^\sigma A(x, f(x)).
\end{equation}
This formulation of the axiom is strictly stronger than the Axiom of Dependent Choices which usually treated in the context of realizability interpretations for Analysis. One can also consider the axioms of Figure~\ref{fig:additional} to be part of \PAomegaplus. It is known that due to their logical form, the proof interpretation that we are going to employ has no modifying effect on them. On the other hand, they can indicate in which sense the realizability model supports extensionality.

\PAomegaplus\ is a first-order, predicate logic which is multi-sorted, that is, has variables and quantifiers that range over the types of System T. The system is ``minimal'' in the sense of Schwichtenberg \cite{SchwichtenbergW2012} -- one wants to know that the method works even if we do not have a special treatment of the absurdity symbol $\bot$; one can work with any fixed formula $N$ as if it were $\bot$. 

\PAomegaplus\ has explicit rules for dealing with the existential quantifier (disjunction has not been included for the sake of simplicity). This, together with the special \textsc{Shift} rule, justifies the plus superscript \textsuperscript{+} in the name. The rule \textsc{Shift} is a general form of the more usual double-negation elimination rule, restricted to $\Sigma_2$-formulas, that is precisely suitable for a simple proof of the Soundness Theorem. It has previously been used by Nakata and the author in a semi-classical logic context \cite{Ilik2010,IlikN2014}.

The class of \emph{$\Sigma_2$-formulas} consists of formulas $S$ of the following form,
\[
S ::= N ~|~ \exists x^\bn N ~|~ N\to S ~|~ N\wedge S ~|~ S\wedge N,
\]
where $N$ stands for so called \emph{computationally irrelevant} formulas \cite{Seisenberger}, defined inductively by
\[
N ::= P ~|~ N \wedge N ~|~ \forall x^\tau N ~|~ A \to N,
\]
where $P$ stands for prime formulas (predicates) and $A$ has no restrictions on the form. 

\begin{figure}
  \centering
  \begin{tabular}{cc}
    ~ & ~ \\
    \axc{$~$}
    \uic{$A,\Gamma\vdash A$}{\textsc{Ax}}
    \DisplayProof
    & 
    \axc{$\Gamma\vdash A$}
    \uic{$B,\Gamma\vdash A$}{\textsc{Wkn}}
    \DisplayProof
    \\
    ~ & ~ \\
    \axc{$A,\Gamma\vdash B$}
    \uic{$\Gamma\vdash A\to B$}{$\to_{\textsc{I}}$}
    \DisplayProof
    &
    \axc{$\Gamma\vdash A\to B$}
    \axc{$\Gamma\vdash A$}
    \bic{$\Gamma\vdash B$}{$\to_{\textsc{E}}$}
    \DisplayProof
    \\
    ~ & ~ \\
    \axc{$\Gamma\vdash A\wedge B$}
    \uic{$\Gamma\vdash A$}{$\wedge^1_{\textsc{E}}$}
    \DisplayProof
    & 
    \axc{$\Gamma\vdash A\wedge B$}
    \uic{$\Gamma\vdash B$}{$\wedge^2_{\textsc{E}}$}
    \DisplayProof
    \\
    ~ & ~ \\
    \axc{$\Gamma\vdash A$}
    \axc{$\Gamma\vdash B$}
    \bic{$\Gamma\vdash A\wedge B$}{$\wedge_{\textsc{I}}$}
    \DisplayProof
    &
    \axc{$\Gamma\vdash A(r^\tau)$}
    \uic{$\Gamma\vdash \exists x^\tau A(x)$}{$\exists_{\textsc{I}}$}
    \DisplayProof
    \\
    ~ & ~ \\
    \multicolumn{2}{c}{
      \axc{$\Gamma\vdash \exists x^\tau A(x)$}
      \axc{$\Gamma\vdash \forall x^\tau (A(x) \to B)$}
      \axc{$x\not\in\FV(B)$}
      \tic{$\Gamma\vdash B$}{$\exists_{\textsc{E}}$}
      \DisplayProof}
    \\
    ~ & ~ \\
    \axc{$\Gamma\vdash A(x^\tau)$}
    \axc{$x\not\in\FV(\Gamma)$}
    \bic{$\Gamma\vdash \forall x^\tau A(x)$}{$\forall_{\textsc{I}}$}
    \DisplayProof
    &
    \axc{$\Gamma\vdash \forall x^\tau A(x)$}
    \uic{$\Gamma\vdash A(r^\tau)$}{$\forall_{\textsc{E}}$}
    \DisplayProof
    \\
    ~ & ~ \\
    \multicolumn{2}{c}{
      \axc{$\Gamma\vdash A(\zero)$}
      \axc{$\Gamma\vdash \forall x^\bn (A(x) \to A(\suc x))$}
      \bic{$\Gamma\vdash \forall x^\bn A(x)$}{\textsc{Ind}}
      \DisplayProof}
    \\
    ~ & ~ \\
    \multicolumn{2}{c}{
      \axc{$\forall x^\bn(A(x)\to S(x)),\Gamma\vdash S(r)$}
      \uic{$\Gamma\vdash A(r)$}{\textsc{Shift}}
      \DisplayProof\qquad ($A,S\in\Sigma_2$)}
    \\
    ~ & ~ \\
  \end{tabular}
  \caption{A natural deduction system for the theory \PAomegaplus}
  \label{fig:haplus}
\end{figure}

\begin{figure}
  \centering
  \begin{tabular}{cc}
    ~ & ~ \\
    \axc{$\reify{~}{}{\eval{}{}{r}{\Uparrow\Gamma}} = \reify{~}{}{\eval{}{}{s}{\Uparrow\Gamma}} $}
    \uic{$\Gamma\vdash r^\tau \doteq s^\tau$}{\textsc{Refl}}
    \DisplayProof
    &
    \axc{$\Gamma\vdash A(r)$}
    \axc{$\Gamma\vdash r \doteq s$}
    \bic{$\Gamma\vdash A(s)$}{\textsc{Comp}}
    \DisplayProof
    \\
    ~ & ~ \\
    \axc{$\Gamma\vdash (\suc r) \doteq \zero$}
    \uic{$\Gamma\vdash \bot$}{\textsc{Cont}}
    \DisplayProof
    &
    \axc{$\Gamma\vdash \bot$}
    \uic{$\Gamma\vdash A$}{\textsc{Efq}}
    \DisplayProof
    \\
    ~ & ~ \\
  \end{tabular}
  \caption{Additional rules for equality (computationally irrelevant)}
  \label{fig:additional}
\end{figure}

As our interpretation we will use a version of Kreisel's so called modified realizability interpretation \cite{Troelstra1973,Kohlenbach2008}, optimised similarly to the work of Berger, Buchholz, Schwichtenberg, and Seisenberger \cite{BergerBS2002,Seisenberger}. However, the key addition is that the realizing terms will also be computed by normalized instances of the $\shift$ term.

\begin{definition}\label{def:mr} Given a context $\Gamma$ and an interpretation of hypotheses $\rho : |\Gamma|\Vvdash|\Gamma|$, the \emph{modified realizability interpretation} ``$\sr{p}{A}$'' of a formula $A$ by a term $p$ of type $|\Gamma|\rvdash|A|$ of System T is defined by the following formula translation,
  \begin{align*}
    \sr{p}{N} &:= N \quad\quad\quad\text{(for any term $p$ of type $\bn$)}\\
    \sr{p}{N \wedge B} &:= N \wedge (\sr{p}{B})\\
    \sr{p}{A \wedge N} &:= (\sr{p}{A}) \wedge N \\
    \sr{p}{A \wedge B} &:= \srConj{p}{A}{B}{\rho}\\
    \sr{p}{N \to B} &:= N \to (\sr{p}{B})\\
    \sr{p}{A \to B} &:= \srImpl{p}{x}{A}{B}{\rho}\\
    \sr{p}{\forall x^\tau A(x)} &:= \srAll{p}{x}{\tau}{A(x)}{\rho}\\
    \sr{p}{\exists x^\tau N(x)} &:= N(p)\\
    \sr{p}{\exists x^\tau A(x)} &:= \srEx{p}{A}{\rho},
  \end{align*}
  in which $N$ denotes a computationally irrelevant formula, and where the type $|A|$ of the realizing term $p$ is computed as follows:
  \begin{align*}
    |N| &:= \bn \\
    |N\wedge B| &:= |B| \\
    |A\wedge N| &:= |A| \\
    |A\wedge B| &:= |A|* |B| \\
    |N\to B| &:= |B| \\
    |A\to B| &:= |A|\to |B| \\
    |\forall x^\tau A| &:= \tau \to |A|\\
    |\exists x^\tau N| &:= \tau\\
    |\exists x^\tau A| &:= \tau * |A|\\
  \end{align*}
  The map $|\cdot|$ is extended to contexts $\Gamma$ by $|C_1,\ldots,C_n| := |C_1|;\cdots;|C_n|$. 
\end{definition}

Note that $\Sigma_2$-formulae are exactly those ones that are realized by a term of type $\bn$.

Our main result is the following one.

\begin{soundthm}
  If \PAomegaplus+AC proves $C_1,C_2,\ldots,C_n\vdash A,$ and $A$ is computationally relevant, then there exists a term $p$ of System \Tplus\ such that \PAomegaplus\ alone proves that, for every $\rho : |C_1|,|C_2|,\ldots,|C_n|\Vvdash |C_1|,|C_2|,\ldots,|C_n|$,
  \begin{equation*}
  \sr{\reify{}{}{\eval{}{}{\Hyp}{\rho}}}{C_1}, \sr{\reify{}{}{\eval{}{}{\Wkn{\Hyp}}{\rho}}}{C_2}, \ldots, \sr{\reify{}{}{\eval{}{}{\wkn^n{\Hyp}}{\rho}}}{C_n} \vdash
  \sr{\reify{}{}{\eval{}{}{p}{\rho}}}{A}.
  \end{equation*}
\end{soundthm}
\begin{proof}
  The proof is by induction on the derivation of $C_1,C_2,\ldots,C_n\vdash A$ and provides the same realizing terms as Kreisel's modified realizability interpretation of \HAomega+AC. The additional rule of \textsc{Shift} (treated in more detail below) is realized via the $\shift$ term of System \Tplus, nevertheless normalized to System T using $\reify{}{}{\eval{}{}{\cdot}{}}$.

  We will denote by $p$ realizing terms provided by the induction hypothesis for each proof rule. If there are two induction hypotheses, the term corresponding to the second, right premise of the proof rule will be denoted by $q$.

  As a general guide, the elimination rules are enough to prove their own soundness, while the introduction rules and the rules AC, \textsc{Wkn}, \textsc{Ind}, and \textsc{Shift} also need to use the definitional equalities of Proposition~\ref{prop:equalities}. 

  In general, the axiom \ref{eq:acsigmatau} is realized by the term \[\Lam{\Pair{\Lam{\App{\Fst{\Wkn{\Hyp}}}{\Hyp}}}{\Lam{\App{\Snd{\Wkn{\Hyp}}}{\Hyp}}}}.\] When the formula $A(x,y)$ is computationally irrelevant, the realizer is the term $\Lam{\Hyp}$. The proof in both cases is a trivial intuitionistic implication and does not require \AC{\sigma}{\tau} itself. Equations (\ref{eq:wkn})-(\ref{eq:lam}) of Proposition~\ref{prop:equalities} are nevertheless used.

  \textsc{Ax} is realized by $\Hyp$.

  \textsc{Wkn} is realized by $\Wkn{p}$ and verified using equation (\ref{eq:wkn}) and (\ref{eq:hyp}).

  The general case of $\to_{\textsc{I}}$, when $A$ from $A\to B$ is computationally relevant, is realized by $\Lam{p}$ and verified using equations (\ref{eq:lam}), (\ref{eq:hyp}), and (\ref{eq:wkn}). The induction hypothesis needs to be used with the context $\rho := (\eval{}{}{x}{\rho} , \rho)$, where $x$ comes from the unfolding of the mr-definition for implication. The special case, when $A$ is computationally irrelevant, is rather realized by the term $p$ only.

  When $A$ from $A\to B$ is computationally relevant, the case $\to_{\textsc{E}}$ is realized by $\App{p}{q}$. When $A$ is irrelevant, the realizing term is just $p$.

  $\wedge^1_{\textsc{E}}$ is realized using $\Fst{p}$, in general, while in the case where one of the conjuncts of $A \wedge B$ is irrelevant, the realizer is just $p$.

  $\wedge^2_{\textsc{E}}$ is realized using $\Snd{p}$, in general, while in the case where one of the conjuncts of $A \wedge B$ is irrelevant, the realizer is just $p$.

  In general, $\wedge_{\textsc{I}}$ is realized using $\Pair{p}{q}$ and verified via equations (\ref{eq:fst}) and (\ref{eq:snd}). When $A$ from $A \wedge B$ is irrelevant, the realizer is $q$, while when $B$ is irrelevant, the realizer is $p$.

  In general, $\exists_{\textsc{I}}$ is realized by $\Pair{p}{q}$ and verified via equations (\ref{eq:fst}) and (\ref{eq:snd}). When $A(x)$ from $\exists x A(x)$ is computationally irrelevant, then the realizer is $r$, the witnessing term.

  $\exists_{\textsc{E}}$ is realized by $\App{\App{q}{\Fst{p}}}{\Snd{p}}$, in general. When $A(x)$ from $\exists x A(x)$ is computationally irrelevant, then the realizer is $\App{q}{p}$.

  $\forall_{\textsc{I}}$ is realized by $\Lam{p}$. For verification, it is necessary to apply equation (\ref{eq:lam}) and to use the induction hypothesis with context $\rho := (\eval{}{}{x}{\rho} , \rho)$.

  $\forall_{\textsc{E}}$ is realized by $\App{p}{r}$.

  \textsc{Ind} is realized by $\Lam{\Rec{\Hyp}{p}{q}}$ and using equations (\ref{eq:zero}) and (\ref{eq:suc}).

  \textsc{Shift} is realized by $\shift{p}$ (normalized to System T using $\reify{}{}{\eval{}{}{-}{}}$). The goal is to prove 
  \begin{equation*}
  \sr{\reify{}{}{\eval{}{}{\Hyp}{\rho}}}{C_1}, \ldots, \sr{\reify{}{}{\eval{}{}{\wkn^n{\Hyp}}{\rho}}}{C_n} \vdash
  \sr{\reify{}{}{\eval{}{}{\Shift{p}}{\rho}}}{A(r)}.
  \end{equation*}
  Using equation (\ref{eq:shift:N}), we obtain $\phi$ and the goal becomes
  \begin{equation*}
  \sr{\reify{}{}{\eval{}{}{\Hyp}{\rho}}}{C_1}, \ldots, \sr{\reify{}{}{\eval{}{}{\wkn^n{\Hyp}}{\rho}}}{C_n} \vdash
  \sr{\reify{}{}{\eval{}{}{p}{\phi,\rho}}}{A(r)}.
  \end{equation*}
  We can now use the induction hypothesis with $\rho:=(\phi,\rho)$,
  \begin{multline*}
  \sr{\reify{}{}{\eval{}{}{\Hyp}{\phi,\rho}}}{\forall x^\bn(A(x)\to S(x))}, \sr{\reify{}{}{\eval{}{}{\wkn\Hyp}{\phi,\rho}}}{C_1}, \ldots, \\\sr{\reify{}{}{\eval{}{}{\wkn^{n+1}{\Hyp}}{\phi,\rho}}}{C_n} \vdash
  \sr{\reify{}{}{\eval{}{}{p}{\phi,\rho}}}{S(r)}.
  \end{multline*}
  Thanks to equation (\ref{eq:wkn}), the induction hypothesis becomes
  \begin{multline*}
  \sr{\reify{}{}{\eval{}{}{\Hyp}{\phi,\rho}}}{\forall x^\bn(A(x)\to S(x))}, \sr{\reify{}{}{\eval{}{}{\Hyp}{\rho}}}{C_1}, \ldots, \\\sr{\reify{}{}{\eval{}{}{\wkn^{n}{\Hyp}}{\rho}}}{C_n} \vdash
  \sr{\reify{}{}{\eval{}{}{p}{\phi,\rho}}}{S(r)}.
  \end{multline*}
  Finally, thanks to equation (\ref{eq:hyp:cont}), we can finish the proof by applying the \textsc{Shift} rule for:
  \begin{align*}
    S'(x,y) &:= \sr{\reify{}{}{\eval{}{}{y}{\phi,\rho}}}{S(x)}\\
    A'(x,y) &:= \sr{\reify{}{}{\eval{}{}{y}{\phi,\rho}}}{A(x)}.
  \end{align*}





\end{proof}

\begin{remark}
  The \textsc{Shift} case in the proof of the Soundness Theorem only uses the case where $\shift$ is reified at type $\bn$. This use does not exhaust the possibilities of the realizability model. For example, one can prove the soundness of the \textsc{Shift} rule for $A = A_1\wedge A_2$ or $A = \exists z^\bn A_2(z)$, when $A_i\in\Sigma_2$, by using the equations
  \begin{align*}
    \reify{}{\tau*\sigma}{\eval{}{}{p}{\rho}} &= \Pair{\reify{}{\tau}{\eval{}{}{\Fst{p}}{\rho}}}{\reify{~}{\sigma}{\eval{}{}{\Snd{p}}{\rho}}} \\
    \reify{}{\bn*\sigma}{\eval{}{}{\Fst{\Shift{p}}}{\rho}} &= \reify{~}{\bn}{\eval{}{}{p}{\phi_1,\rho}}\\
    \reify{}{\tau*\bn}{\eval{}{}{\Snd{\Shift{p}}}{\rho}} &= \reify{~}{\bn}{\eval{}{}{p}{\phi_2,\rho}}\\
    \reify{}{\bn}{\eval{}{}{\App{\App{\hyp}{x}}{y}}{\phi_1,\rho}} &= \reify{~}{\bn}{\eval{}{}{\Fst{y}}{\phi_1,\rho}}\\
    \reify{}{\bn}{\eval{}{}{\App{\App{\hyp}{x}}{y}}{\phi_2,\rho}} &= \reify{~}{\bn}{\eval{}{}{\Snd{y}}{\phi_2,\rho}}
  \end{align*}
  where
  \begin{align*}
    x & : \bn\to\bn*\bn\to\bn;\gamma \vdash \bn\\
    y & : \bn\to\bn*\bn\to\bn;\gamma \vdash \bn*\bn\\
    \phi_1 &:= \eta (\ge_2 \nu \mapsto \eta (\ge_3 \alpha \mapsto \eta (\alpha \gerefl (\ge_4 \gamma \mapsto (\proj_1 \gamma)))))\\
    \phi_2 &:= \eta (\ge_2 \nu \mapsto \eta (\ge_3 \alpha \mapsto \eta (\alpha \gerefl (\ge_4 \gamma \mapsto (\proj_2 \gamma))))).
  \end{align*}
  The induction hypothesis needs to be used twice, once for $\phi_1$ and once for $\phi_2$.

  Similar equations hold for function types:
  \begin{gather*}
    \reify{\gamma}{\bn}{\eval{}{}{\App{\Shift{p}}{z}}{\rho}} = \reify{\bn\to(\bn\to\bn)\to\bn;\gamma}{\bn}{\eval{}{}{p}{\phi_3, \rho}}\\
    \phi_3 := \eta (\ge_2 \nu \mapsto \eta (\ge_3 \alpha \mapsto \eta (\alpha \gerefl (\ge_4 \gamma \mapsto \mu(\gamma \gerefl \eval{}{}{z}{\Vvdashge{\ge_4\cdot\ge_3\cdot\ge_2}{\rho}})))))
  \end{gather*}
  Nevertheless, it does not appear to be possible to prove the soundness of the \textsc{Shift} rule even for the case $|A|=\bn\to\bn$ in general, since it is well known that there are already classically true $\Sigma^0_3$-formulas which do not have a recursive realizer. Still, it may be the case that the realizability model can be used to give a sound computational interpretation of \emph{particular} $\Sigma^0_3$ (or more complex) formulas.
\end{remark}

\begin{corollary}\label{cor:gamma1}
  The $\Sigma_2$-fragment of classical Analysis satisfies the Existence Property,
  \begin{quote}
    Given a derivation of $\Gamma\vdash\exists x^\tau A(x)$, there exists a term $p$ of type $\tau$ of System T such that $\Gamma\vdash A(p)$.
  \end{quote}
  and, consequently, the Weak Church's Rule,
  \begin{quote}
    Given a (closed) derivation of $\emptyset\vdash\forall x^\bn\exists y^\bn A(x,y)$, there exists a total recursive function $\mathbf{f} : \bn \to \bn$ such that, for all $\mathbf{n}\in\bn$, we have that $\emptyset\vdash A(\overline{\mathbf{n}},\overline{\mathbf{f n}})$, where $\overline{\mathbf{m}}$ denotes the term $\underbrace{\suc\cdots\suc}_{\mathbf{m} \text{ times }}\zero$.
  \end{quote}
\end{corollary}
\begin{proof}
The proof method is not new (see Corollary~5.24 of \cite{Kohlenbach2008} and paragraph~1.11.7 of \cite{Troelstra1973}). If the formula $A$ is of the class $\Gamma_1$ \cite{Troelstra1973,Kohlenbach2008},
\[
\Gamma_1 \ni G ::= N ~|~ G\wedge G ~|~ \forall x G ~|~ \exists x G ~|~ S \to G,
\]
(where $S$ is a $\Sigma_2$-formula and $N$ is computationally irrelevant) then already intuitionistic logic shows that $\vdash (\sr{p}{A}) \to A$. If $A$ is outside this class, one first needs to define the ``with truth'' variant of modified realizability in which one replaces the clause for implication of mr-interpretation of Definition~\ref{def:mr} by
\[
\srt{p}{A \to B} := \srtImpl{p}{x}{A}{B}{\rho} \wedge (A \to B).
\]
The Soundness Theorem is provable for mrt with the same realizing terms, but now we also have, for any formula $A$, $(\srt{p}{A}) \to A$. This directly implies the Existence Property.

For the special case when $\Gamma=\emptyset$ and $\tau=\bn$, we get the Numerical Existence Property:
\begin{quote}
    Given a (closed) derivation of $\emptyset\vdash\exists y^\bn A(y)$, there exists $\mathbf{n} : \bn$ such that $\emptyset\vdash A(\overline{\mathbf{n}})$.
\end{quote}
This follows from the fact that a closed derivation has a realizer that is a closed term (does not have non-bound $\hyp$ subterms). Since the realizer is necessarily in normal form, and since it is not neutral (all neutral terms have at least one non-bound occurrence of $\hyp$), then the realizer must be of the required form $\underbrace{\suc\cdots\suc}_{\mathbf{n} \text{ times }}\zero$.

To show the Weak Church's Rule, we use elementary Recursion Theory. Like all theories over countable languages, \PAomegaplus+AC is recursively axiomatizable, that is, there exists a recursive predicate $\Proof(k,l)$ formalizing the fact that $k\in\bn$ is a code for a derivation of the formula coded by $l\in\bn$. 

Let $g(n) = \min_m \Proof(j_1 m, \ulcorner A(\overline{n}, \overline{j_2 m})\urcorner)$, where $j_1$ and $j_2$ are the projections of some surjective pairing function. As defined, $g$ is a partial recursive function.

Now, given $\emptyset \vdash \forall x^\bn\exists y^\bn A(x,y)$ and $n\in\bn$, we obtain $\emptyset \vdash \exists y^\bn A(\overline{n}, y)$, and by the Numerical Existence Property we obtain $m\in\bn$ such that $\vdash A(\overline{n},\overline{m})$. We proved that, for every $n$, there exists $m$ such that $\emptyset \vdash A(\overline{n},\overline{m})$ which shows that the function $g$ is \emph{total} recursive. We may now take $f(n) := j_2(g(n))$ and by definition we have that, for any $n$, $\emptyset \vdash A(\overline{n},\overline{f(n)})$.  
\end{proof}

Note that the class $\Sigma_2$ includes the following schemata,
\begin{equation}
  \label{eq:mp}
  \tag{MP}
  \neg_N \neg_N \exists x^\bn M(x) \to \exists x^\bn M(x),
\end{equation}
\begin{equation}
  \label{eq:dns0}
  \tag{DNS}
  \forall x^\bn\neg_N\neg_N A(x) \to \neg_N\neg_N\forall x^\bn A(x),
\end{equation}
where $M, N$ denote computationally irrelevant formulas and $\neg_N A$ denotes negation in minimal logic, that is $A\to N$ for a fixed $N$.

The Existence Property implies that principles like \ref{eq:mp} and \ref{eq:dns0} can justly be considered as constructive even in presence of induction and the full Axiom of Choice, partly extending previous works \cite{Troelstra1973,Gabbay1972,Seldin1986,HerbelinMP,Ilik2010}. Similar conclusions follow from the work of Rand Moschovakis that uses a version of Kleene's general-recursive realizability \cite{Moschovakis2002}.

Weak Church's Rule seems to justify why, even constructively, \ref{eq:ct0} deserves the name ``the false Church's Thesis'' \cite{Lombardi}.

\section*{Acknowledgments}
This work was funded by Kurt Gödel Research Prize Fellowship 2011 and ERC Advanced Grant ProofCert. I would also like to thank Mart\'in Escard\'o, Jaime Gaspar, Keiko Nakata, and Dirk Pattinson for comments on an earlier version of this paper, and Dale Miller for providing scientific liberty.

\bibliographystyle{plain}
\bibliography{shift-analysis}

\begin{thebibliography}{10}

\bibitem{Berger2004}
Ulrich Berger.
\newblock A computational interpretation of open induction.
\newblock In {\em Proceedings of the 19th Annual IEEE Symposium on Logic in
  Computer Science (LICS'04)}, pages 326--334. IEEE Computer Society, 2004.

\bibitem{BergerBS2002}
Ulrich Berger, Wilfried Buchholz, and Helmut Schwichtenberg.
\newblock Refined program extraction from classical proofs.
\newblock {\em Annals of Pure and Applied Logic}, 114:3--25, 2002.

\bibitem{BergerO2005}
Ulrich Berger and Paulo Oliva.
\newblock Modified bar recursion and classical dependent choice.
\newblock {\em Springer-Verlag Lecture Notes in Logic}, 20:89--107, 2005.

\bibitem{BergerO2006}
Ulrich Berger and Paulo Oliva.
\newblock Modified bar recursion.
\newblock {\em Mathematical Structures in Computer Science}, 16:163--183, 2006.

\bibitem{BergerS1991}
Ulrich Berger and Helmut Schwichtenberg.
\newblock An inverse of the evaluation functional for typed lambda-calculus.
\newblock In {\em LICS}, pages 203--211. IEEE Computer Society, 1991.

\bibitem{Danvy1996}
Olivier Danvy.
\newblock Type-directed partial evaluation.
\newblock In {\em POPL}, pages 242--257, 1996.

\bibitem{DanvyF1990}
Olivier Danvy and Andrzej Filinski.
\newblock Abstracting control.
\newblock In {\em LISP and Functional Programming}, pages 151--160, 1990.

\bibitem{Gabbay1972}
Dov~M. Gabbay.
\newblock Applications of trees to intermediate logics.
\newblock {\em The Journal of Symbolic Logic}, 37:135--138, 1972.

\bibitem{GodelPostscriptSpector}
Kurt G{\"o}del.
\newblock {\em Collected works. Publications 1938--1974}, volume~II, chapter
  Postscript to {S}pector 1962, page 253.
\newblock The Clarendon Press Oxford University Press, New York, 1962.

\bibitem{Godel1958}
Kurt G{\"o}del.
\newblock {\em Collected works. Publications 1938--1974}, volume~II, chapter On
  a hitherto unutilized extension of the finitary standpoint, pages 241--251.
\newblock The Clarendon Press Oxford University Press, New York, 1990.
\newblock (English translation of the original 1958 article).

\bibitem{Godel1941}
Kurt G{\"o}del.
\newblock {\em Collected works. Unpublished essays and lectures}, volume III,
  chapter In what sense is intuitionistic logic constructive, pages 189--200.
\newblock The Clarendon Press Oxford University Press, 1995.
\newblock (early lecture on the Dialectica interpretation from 1941).

\bibitem{HerbelinMP}
Hugo Herbelin.
\newblock An intuitionistic logic that proves {M}arkov's principle.
\newblock In {\em Proceedings of the 25th Annual IEEE Symposium on Logic in
  Computer Science, LICS 2010, 11-14 July 2010, Edinburgh, United Kingdom},
  pages 50--56. IEEE Computer Society, 2010.

\bibitem{Herbelin2012}
Hugo Herbelin.
\newblock A constructive proof of dependent choice, compatible with classical
  logic.
\newblock In {\em Proceedings of the 27th Annual ACM/IEEE Symposium on Logic in
  Computer Science, LICS 2012, 25-28 June 2012, Dubrovnik, Croatia}, pages 365
  -- 374. IEEE Computer Society, 2012.

\bibitem{Howard1968}
William~Alvin Howard.
\newblock Functional interpretation of bar induction by bar recursion.
\newblock {\em Composition Mathematica}, 20:107--124, 1968.

\bibitem{IlikThesis}
Danko Ilik.
\newblock {\em Constructive Completeness Proofs and Delimited Control}.
\newblock PhD thesis, École Polytechnique, October 2010.

\bibitem{Ilik2010}
Danko Ilik.
\newblock Delimited control operators prove double-negation shift.
\newblock {\em Annals of Pure and Applied Logic}, 163(11):1549 -- 1559, 2012.

\bibitem{Ilik2013}
Danko Ilik.
\newblock Continuation-passing style models complete for intuitionistic logic.
\newblock {\em Annals of Pure and Applied Logic}, 164(6):651 -- 662, 2013.

\bibitem{code}
Danko Ilik.
\newblock Formal proof of normalization of {S}ystem {T+} in {A}gda.
\newblock Available at
  \url{http://www.lix.polytechnique.fr/~danko/shift-analysis.zip}, 2014.

\bibitem{IlikN2014}
Danko Ilik and Keiko Nakata.
\newblock A direct version of {V}eldman's proof of open induction on {C}antor
  space via delimited control operators.
\newblock {\em Leibniz International Proceedings in Informatics (LIPIcs)},
  26:188--201, 2014.

\bibitem{Kohlenbach2008}
Ulrich Kohlenbach.
\newblock {\em Applied Proof Theory: Proof Interpretations and Their Use in
  Mathematics}.
\newblock Springer Monographs in Mathematics. Springer-Verlag, Berlin,
  Heidelberg, 2008.

\bibitem{Kreisel1959}
Georg Kreisel.
\newblock Interpretation of analysis by means of constructive functionals of
  finite types.
\newblock In Arend Heyting, editor, {\em Constructivity in Mathematics,
  Proceedings of the colloqium held at Amsterdam, 1957}, Studies in Logic and
  The Foundations of Mathematics, pages 101--127. North-Holland Publishing
  Company Amsterdam, 1959.

\bibitem{KreiselBR1976}
Georg Kreisel.
\newblock Review of the paper ``{T}he model {G} of the theory {BR}'' by
  {E}rsov.
\newblock {\em Zentralblatt für Mathematik und ihre Grenzgebiete}, 312, 1976.

\bibitem{Krivine2003}
Jean-Louis Krivine.
\newblock Dependent choice, ‘quote’ and the clock.
\newblock {\em Theoretical Computer Science}, 308(1–3):259 -- 276, 2003.

\bibitem{Kuroda1951}
Sigekatu Kuroda.
\newblock Intuitionistische untersuchungen der formalistischer logik.
\newblock {\em Nagoya Mathematical Journal}, 3:35--47, 1951.

\bibitem{Lombardi}
Henri Lombardi and Claude Quitté.
\newblock {\em Algèbre commutative -- Méthodes constructives}.
\newblock Calvage \& Mounet, Paris, 2011.

\bibitem{Moschovakis2002}
Joan~Rand Moschovakis.
\newblock Analyzing realizability by {T}roelstra's methods.
\newblock {\em Annals of Pure and Applied Logic}, 114:203--225, 2002.

\bibitem{Schwichtenberg1979}
Helmut Schwichtenberg.
\newblock On bar recursion of types 0 and 1.
\newblock {\em The Journal of Symbolic Logic}, 44(3), 1979.

\bibitem{SchwichtenbergW2012}
Helmut Schwichtenberg and Stanley~S. Wainer.
\newblock {\em Proofs and Computations}.
\newblock Perspectives in Logic. Cambridge University Press, 2012.

\bibitem{Seisenberger}
Monika Seisenberger.
\newblock Program from proofs using classical dependent choice.
\newblock {\em Annals of Pure and Applied Logic}, 153:97--110, 2008.

\bibitem{Seldin1986}
Jonathan~P. Seldin.
\newblock On the proof theory of the intermediate logic {MH}.
\newblock {\em The Journal of Symbolic Logic}, 51(3):626--647, 1986.

\bibitem{Spector1962}
Clifford Spector.
\newblock Provably recursive functionals of analysis: a consistency proof of
  analysis by an extension of principles formulated in current intuitionistic
  mathematics.
\newblock In {\em Proc. {S}ympos. {P}ure {M}ath., {V}ol. {V}}, pages 1--27.
  American Mathematical Society, Providence, R.I., 1962.

\bibitem{BerardiBC1998}
Marc~Bezem Stefano~Berardi and Thierry Coquand.
\newblock On the computational content of the axiom of choice.
\newblock {\em The Journal of Symbolic Logic}, 63(2):600--622, 1998.

\bibitem{Troelstra1973}
Anne~S. Troelstra, editor.
\newblock {\em Metamathematical Investigations of Intuitionistic Arithmetic and
  Analysis}.
\newblock Number 344 in Lecture Notes in Mathematics. Springer-Verlag, Berlin,
  Heidelberg, New York, 1973.

\bibitem{hottbook}
The {Univalent Foundations Program}.
\newblock {\em Homotopy Type Theory: Univalent Foundations of Mathematics}.
\newblock Institute for Advanced Study, 2013.

\end{thebibliography}

\end{document}